# Communication of Prospective Teachers with Students in Mathematics Learning at Senior High School (SMA)

Dian Kurniawan, Ipung Yuwono, Edy Bambang Irawan, Hery Susanto, Subanji, Susiswo.

Mathematics Education Department, Faculty of Post Graduate, State University of Malang, Indonesia.
e-mail: dian.kurniawan@unsil.ac.id

**ABSTRACT**

This study aims to describe the communication of prospective teachers in learning mathematics. The research was conducted in SMAN 1 Madura. The subject of this research is PLP student from Mathematics Education Study Program of Madura University. The prospective teacher's math learning is using the concept of function limits and the derivative of functions in problem-solving (Calculus). This research is descriptive qualitative research. Instruments in the study were interviews, field notes, and video recordings. The knowledge of prospective mathematics teacher in field practice involves communication in factual, conceptual, procedural, and metacognitive knowledge in facts, concepts, principles, procedures, and reasoning. Teacher communication plays a role in learning mathematics, especially the ability to use formulas, rules, and methods simultaneously correctly to make accurate calculations. Prospective teachers are able to provide guidance to their students when implementing mathematics lessons and when working on the blackboard, directing discussions and frequently asked questions, using language appropriate to the level of students' understanding so that the symbols, procedures, and strategies described can be well understood by the students. Factors that are considered important to be considered in the learning of mathematics in an effort to develop the skills of prospective teachers so that students are accustomed to facing similar problems, and apply a mathematical procedure in a newly encountered context. The result of this research is expected for prospective mathematics teacher to master communication in learning mathematics. The communication of prospective teachers in mathematics learning can be applied in other wider aspects.

1. **INTRODUCTION**

Knowledge of teachers in mathematics subjects needs to be improved through research. The research was conducted in order to develop the knowledge of mathematics teachers on the aspects of teaching and learning. Research is carried out in support of changes in the conception of teachers, teachers' beliefs about their profession, knowledge of their students, curriculum or teaching strategies used, subjects and information submitted by teachers so that teachers need to prepare the material beforehand so that the delivery of accepted teaching facts is acceptable.

(Ball & McDiarmid, 1990) points out: "Basic notes and tricks, and definitions do not find meaningful understanding.", Basic and secondary education has difficulty remembering ideas and procedures so that many do not understand the conceptual understanding of the mathematics they have learned and when delivered. Therefore, teachers should master the aspects of mathematical knowledge well. Aspects The mathematical knowledge that must be mastered involves more detailed conceptual knowledge, the mastery of definitions, conventions (changes), and procedures, finding patterns, conjectures, justifying statements and proving solutions, then seeking generalizations.

The teacher is a professional in his position as an educator. Practical experience in learning required prospective teachers in improving their ability as a prospective educator. (Canrinus, Helms-Lorenz, Beijaard, Buitink, & Hofman, 2011) describe the professional identity of teachers in this regard relating to sustained interpretation and interaction in the context of learning, teachers, commitment to work, changes in confidence and motivation. Learning practices provide useful opportunities to support the growth and development of math or prospective teachers (Amador, 2017).

The knowledge of mathematics prospective teacher is limited mainly because it is based on their experience as students. A prospective teacher must be competent in order to foster students' beliefs so they can receive mathematical material that they learn well. Student confidence will grow if communication between teachers and students in depth in the process of mathematics education is well established.

Growth and development of mathematics skills of prospective teachers have an influence in learning. (Hall, 2014) explains that teacher influences play a role in learning, such as: (1) tend to be energetic, passionate, caring and flexible; (2) sensitive to individual needs and student motivation; (3) enthusiastic and passionate about their teaching subject; and (4) concerned with the value of each student as a person. (Olteanu, 2015) illustrates five teacher approaches to grow and develop as teachers: (1) teachers are meant to develop abilities in the subject matter he taught; (2) master learning strategy; (3) teachers are skilled in delivering the subject matter; (4) effective in teaching; and (5) effective in facilitating student learning.

Learning is done in order to facilitate students to build new knowledge from previous experience and knowledge. Knowledge is built on regular learning, (Sriraman, n.d.) argues that "when students are given what they are given an opportunity to learn.", When students are given the opportunity to learn through their learning activities will receive the transfer of knowledge provided by the teacher, in which case the teacher facilitates the students to gain new knowledge through the learning activities undertaken, and when repetition of the information is delivered either through the assignment, training, and subsequent learning then they will apply prior knowledge. (Sub-county, 2016) affirms that learning is an active process that involves discussion and allows students to reach their own conclusion goals.

Lessons learned should receive feedback from students. This is done through question and answer activities in class and assignment or training, usually with the help of textbooks or textbooks that have been used as a reference in learning. In this feedback activity, communication will occur between teachers and students in the form of oral and written communication. Communication in the form of oral discourse and student writing should not be underestimated. (Pourdavood & Wachira, 2015) explains that the discourse of oral and written communications in the classroom meets three broad and interrelated issues of learning, teaching and assessment purposes.

Learning objectives, Teaching, and Assessment in the form of Oral and Observational communication observed involve observations between the Teacher and students, fellow students, and individual students themselves. (Sür & Delice, 2016) suggest a mathematical communication process evaluated according to individual relationships with each other in a classroom environment, observed observational data and consequently there are three types emerging;

1. Communication process conducted by students and teachers (Teacher Communication to Students).
2. Communication process of students appears each other (Student to student communication).
3. The process of communication that students do with themselves (Internal Communication students).

Evaluation is done depending on the content of the teaching, among others the grouping of objects, which form the mathematical symbolic structure when mathematical communication is used, the communication process occurring in oral statements, there is a situation that oral statements support mathematical objects. Mathematical objects are used in written communication, verbal descriptions of verbal communication. The results of observation of researchers when guiding Field Activity Practices (PLP) prospective Teachers at the University of Siliwangi Tasikmalaya showed one mathematical object in the mathematics learning in the form of skills of the procedures in learning mathematics still have shortcomings especially this is because Candidate Teachers often use algorithm procedures incomplete meaning steps - the process of solving the problem often uses a less raw way quickly. This is in accordance with the National Assessment of Educational Progress (Ball & McDiarmid, 1990) that students are able to perform routine arithmetic calculations, but many have difficulty with reasonably complex procedures and reasoning.

(Level & Teachers, 2003) describes communication requiring greater attention to mathematics education especially in written and oral communication. Students at all grade levels are able to express written and oral communication in order to organize and consolidate their mathematical thinking through communication; communicating their mathematical thoughts coherently and clearly to peers, teachers, and others; analyze and evaluate the mathematical thinking and strategies of others; and use mathematical language to express mathematical ideas appropriately. Communication is an integral part of the classroom and school processes, and the quality of communication affects the quality of teaching and learning of mathematics.(Olteanu, 2015).

The Competency Standards of Primary and Secondary Education graduates in the Regulation of the Minister of Education and Culture No. 20 of 2016 explain that teachers should have factual, conceptual, procedural, and metacognitive knowledge at the technical, specific, detailed and complex levels regarding science, technology, art, culture, and humanities. Procedural Knowledge: Technical and Specific Knowledge, Details and Complex with respect to science, technology, art, and culture related to society and surrounding natural environment, nation, country, regional, and international. Conceptual Knowledge: Terminology / terminology and classification, categories, principles, generalizations, theories, models, and structures used in connection with technical and specific knowledge, detailed and complex with respect to science, technology, art and culture related to society and the natural environment around, nation, state, regional, and international. Procedural Knowledge is Knowledge of how to do things or activities related to technical knowledge, specific, algorithms, methods, and criteria to determine appropriate procedures with regard to science, technology, art and culture, related to society and the natural environment, , country, regional, and international.

("(Mathematics Education Library 10) Alan J. Bishop, Stieg Mellin-Olsen (auth.), Alan J. Bishop, Stieg Mellin-Olsen, Joop Van Dormolen (eds.)-Mathematic.pdf," n.d.) describes the teaching of mathematics, "teaching of mathematics, many different things have been systematically scrutinized, and for the many different purposes, and as the teaching of mathematics is a very complicated matter. ". Mathematics teaching has many different things and is systematic to be studied, with a variety of purposes, as a very complicated problem, and always has a lot to do. Therefore, in this case, the researchers took the title "Communication Prospective Teachers with Students in Mathematics Learning in Senior High School (SMA)".

## 2. METHOD

This research is a qualitative descriptive research. (Creswell, 2009) Data collection techniques conducted in the form of field notes, interviews, documentation, and video recording.

## 3. RESULTS AND DISCUSSION

(Wood, 2012) suggests communication in mathematics learning, one of which is oral communication, presentation of material, that is how to speak in semi-formal form and processing Mathematics idea with a combination of listening and speaking. (Education, 2015) revealed that most teachers believe that success in mathematics means the ability to use formulas, rules, and methods simultaneously correctly to make accurate calculations.

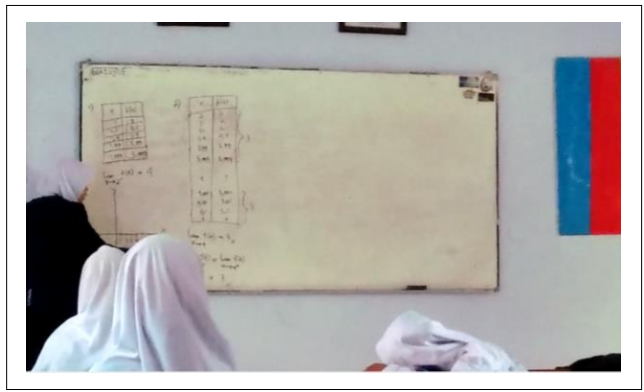

Figure 3.1. Prospective Teacher learning activities.

Prospective teachers appear to have good skills in using formulas, rules, and methods simultaneously correctly to make accurate calculations with Competency Standards using the concept of function limits and function derivatives in problem-solving (Calculus). Prospective teachers communicate in writing while working on questions on the board and explain graphics or drawings.

Prospective teachers are able to coach students while performing mathematics lessons and while working on problems on the board.

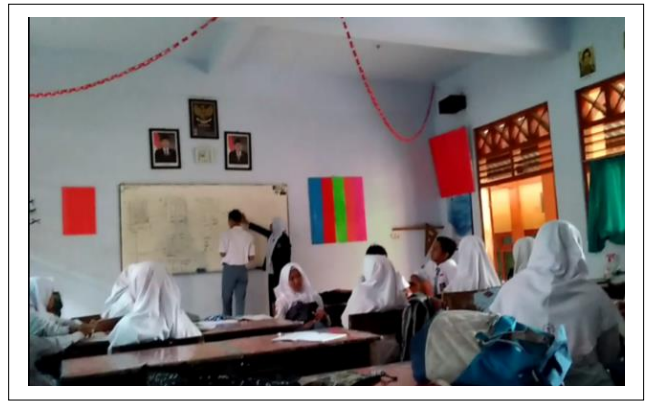

Figure 3.2. Prospective teacher tutoring in working on a problem.

(Guveli, 2015) argues that teachers are people who develop, direct, motivate themselves, develop and practice activities, examine, make students ask and think, listen, work together and evaluate. Coaching prospective teachers in learning is a means of developing student skills.

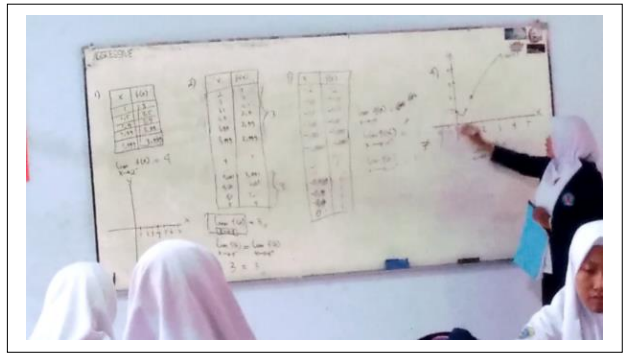

Figure 3.3. An activity describes prospective teachers.

Prospective teachers use good language when explaining the subject matter so as to adjust to the level of student understanding so that the symbols, procedures, and strategies described can be understood by the students. This shows the teachers are skilled in delivering the subject matter.

Communication of prospective teachers is able to stimulate students to be active in the discussion so that good cooperation between students in the discussion process. This is a skill that must be mastered by prospective teachers in implementing learning effectively.

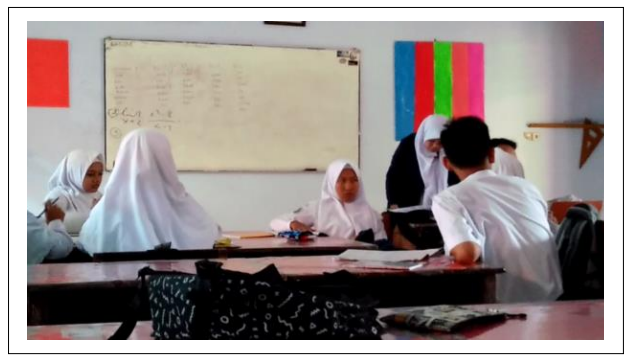

Figure 3.4. Discussion activities.

Prospective teachers are able to direct the question and answer activities with their students so as to establish good communication between prospective teachers and students. This helps the students accept the concept of the subject matter better. Prospective teachers are effective in facilitating student learning.

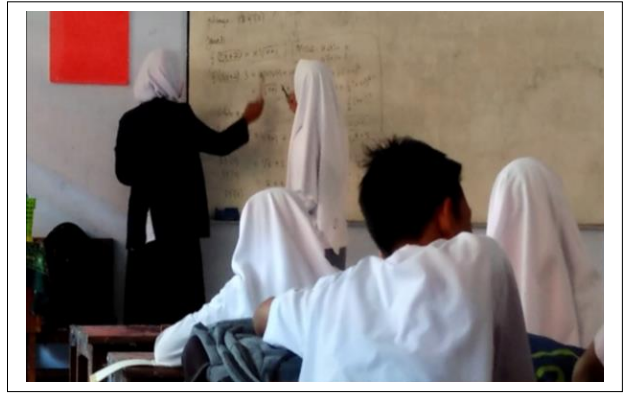

Figure 3.5. Teacher question and answer activities students

## 4. CONCLUSIONS AND RECOMMENDATIONS

**Conclusion**

1. Teacher communication plays a role in learning mathematics especially the ability to use formulas, rules, and methods simultaneously correctly to make accurate calculations.
2. Prospective teachers are able to conduct guidance to their students when implementing mathematics learning and when working on the problem on the blackboard.
3. The language used by the prospective teacher can be adjusted to the level of student's understanding so that the symbols, procedures, and strategies described can be well understood by the students.
4. Prospective teachers are able to direct discussion and question and answer activities well.
5. Factors that are considered important to be considered in the learning of mathematics in an effort to develop the skills of prospective teachers so that students are accustomed to facing similar problems, and apply a mathematical procedure in a newly encountered context.

**Suggestions**

1. Prospective teachers master communication in learning mathematics.
2. Communication of prospective teachers in learning mathematics can be applied in other aspects of a wider